\documentclass[titlepage,11pt]{article}
\usepackage{latexsym}
\usepackage{amsfonts}
\usepackage{amsmath}
\usepackage{tikz}
\usepackage{float}
\usepackage{comment}
\usetikzlibrary{decorations.pathreplacing,calligraphy}
\usetikzlibrary{positioning,arrows.meta}

\newcommand\blackslug{\hbox{\hskip 1pt \vrule width 4pt height 8pt depth 1.5pt
        \hskip 1pt}}
\newcommand\bbox{\hfill \quad \blackslug \bigbreak}

\def\DD{\hbox{-}}
\def\CC{\hbox{-}\cdots\hbox{-}}
\def\LL{,\ldots,}
\def\cupcup{\cup\cdots\cup}

\title{When all directed cycles have length three}
\author{
Paul Seymour\thanks{Supported by AFOSR grant
FA9550-22-1-0234, and NSF grant  DMS-2154169.}\\
Princeton University, Princeton, NJ 08544}

\date{September 17, 2024; revised \today}

\newtheorem{thm}{}[section]

\newcommand{\Proof}{\noindent{\bf Proof.}\ \ }

\begin{document}
\maketitle
\begin{abstract}
We give a construction to build all digraphs with the property that every directed cycle has length three.
\end{abstract}

\section{Introduction}
A {\em digraph} (in this paper) means a directed graph without loops or parallel edges (but we permit directed cyles of length two). A digraph is {\em $\ell$-cyclic} if all its
directed cycles have length $\ell$.
This paper is an exploration of the land of 
3-cyclic digraphs. We will give a construction for all 3-cyclic digraphs, and look at 
some natural subclasses. Some of the easier results work for general values of~$\ell$.

This question arose from trying to extend a theorem of Thomassen~\cite{thomassenacyclic}. This said:
\begin{thm}\label{thomassenthm}
Let $G$ be an acyclic digraph, and let $S=\{a,b\}$ be its set of sources, and $T=\{c,d\}$ its set of sinks.
Suppose that $S\cap T=\emptyset$, and every vertex not in $S\cup T$ has in-degree and out-degree both at least two. Then exactly one of the following holds:
\begin{itemize}
\item there are vertex-disjoint dipaths $P,Q$ of $G$, where $P$ is from $a$ to $c$ and $Q$ is from $b$ to $d$;
\item $G$ can be drawn in a closed disc with $a,b,c,d$ drawn in the boundary in order.
\end{itemize}
\end{thm}
(A shorter proof of a slightly  more general theorem is given in~\cite{acyclic}.) We can think of this as saying that if $G$
is drawn in a rectangle with crossings allowed, with the sources on the north side and sinks on the south side of the rectangle, 
so that $a,b,c,d$ are in order in the boundary of the rectangle, then either
we can rearrange the drawing to eliminated the crossings, without moving the sources or sinks, or there are paths $P,Q$ as in the 
theorem that show this is impossible. Now consider what happens if we glue the north side of the rectangle to the south side, 
making an annulus, and making the identifications $a=d$ and $b=c$. We start with a drawing with crossings in the annulus, and we would 
like to rearrange the drawing (by some continuous deformation, remaining within the annulus at eash stage) to get rid of the crossings; and we
can do so if and only if every directed cycle winds exactly once around the annulus. But there is a ``cut'' across the annulus passing 
through only two vertices, and it is not clear whether this is important. Could the following be true?
\begin{thm}\label{annconj}
{\bf Conjecture:} Let $G$ be a strongly 2-connected digraph drawn (possibly with crossings) in an annulus $\Sigma$,
such that every directed cycle has winding number one around the annulus in the natural sense. Then $G$ can be drawn in $\Sigma$
without crossings, such that every directed cycle has winding number one.
\end{thm}
This would be a pretty extension of Thomassen's theorem, but we cannot prove it.

Drawing in an annulus is the same as drawing 
in the plane
without using the origin; and if we have such a drawing, each edge $e$ of $G$ subtends an angle at the origin (a $w(e)$ fraction of a full 
rotation, say), and for 
each directed cycle $C$, $\sum_{e\in E(C)}w(e)=1$. Let us say a digraph $G$ is {\em weightable} if we can assign a real weight
$w(e)$ to each edge $e$, such that $\sum_{e\in E(C)}w(e)=1$ for each directed cycle $C$. Being weightable is equivalent to being 
drawable (possibly with crossings) in an annulus as in the hypothesis of \ref{annconj}. So, which digraphs are weightable?

In the final section we will show that this is true
if and only if no subgraph of $G$ is a ``weak $k$-double cycle'' with $k\ge 3$ (for the definition, see~\cite{evendigraph}), using
a proof virtually identical to the proof in~\cite{evendigraph}, working with the real field instead of $GF(2)$.
So we can characterize these digraphs 
by excluded
subdigraphs, and hence obtain a co-NP characterization. What about an NP-characterization, can we give a construction for all weightable digraphs? For strongly 2-connected digraphs, \ref{annconj} would give such a construction, 
but that remains open. (Here is one curious fact about them, however: every weightable digraph admits a weighting $w$ that that is $\{0,1\}$-valued. We give the proof in the final section.)

So we look for evidence, simple cases of \ref{annconj} that we can check.
For instance, what if all the numbers $w(e)$  are equal, that is, $G$ is $\ell$-cyclic (for some $\ell$)? Unfortunately, this provides 
no confirmation for \ref{annconj}, because there are no such digraphs: one can show that no $\ell$-cyclic digraph is strongly 2-connected. 
(That is easy: if $G$ has minimum out-degree at least two, construct a directed walk $v_1,v_2\LL $ where for each $i>\ell$, $v_i\ne v_{i-\ell}$. When this repeats a vertex 
we have found a cycle of length different from $\ell$.)  At the moment, that is the extent of progress on \ref{annconj}: we are stuck on 
obtaining a general construction for all weightable digraphs, and even for those that are strongly 2-connected.

But what about $\ell$-cyclic digraphs, can we understand them? Since $\ell$-cyclic digraphs cannot be strongly 2-connected, this suggests that there might be a 
structure theorem for them; can all $\ell$-cyclic digraphs be built from simple ones by some assembly operations? This is so when 
$\ell=2$ (every strongly-connected 2-cyclic digraph is obtained from a tree by replacing each edge by a directed cycle of length two --- thanks to Stephen Bartell for this observation), 
so the first interesting case is when $\ell=3$, the topic of this paper. We will give such a construction when $\ell=3$. In joint work with
Bartell, we tried to do the same when $\ell=4$, but without success so far.

If $G$ is a digraph without directed cycles of length two, let $G^-$ be the graph obtained by forgetting the direction of the edges; that is,
$G^-$ is the graph with vertex set $V(G)$, in which distinct $u,v$ are adjacent if and only if either $v$ is 
adjacent from $u$ in $G$ or vice versa.
We say a non-null digraph $G$ is {\em strongly connected} if for every two vertices $u,v$, there is a directed path from $u$ to $v$;
and {\em weakly connected} if $G^-$ is connected. Similarly, $G$ is {\em weakly $k$-connected} if $G^-$ is $k$-connected. We write
$|G|$ for $|V(G)|$. 

Clearly a digraph $G$ is $\ell$-cyclic if and only if all its strong components are $\ell$-cyclic, so the problem of constructing all $\ell$-cyclic 
digraphs immediately reduces 
to the subproblem of constructing those that are strongly connected. Moreover, if $G_1,G_2$ are $\ell$-cyclic, then we may identify
one vertex of $G_1$ with one vertex of $G_2$, and the result is still $\ell$-cyclic, so it suffices to construct the
$\ell$-cyclic digraphs that are weakly 2-connected. Let us say a digraph $G$ is {\em unbreakable} if it is strongly connected and weakly 2-connected.

Another thing we can do to reduce the question to smaller classes of digraphs is to use the fact that
if $G$ is strongly connected and $\ell$-cyclic, then there is an ordered partition $(A_1\LL A_{\ell})$ of $V(G)$ 
such that 
for every edge $uv$ of $G$ there exists $i\in \{1\LL \ell\}$ such that $u\in A_i$ and $v\in A_{i+1}$ (reading subcripts modulo $\ell$). 
Let us call such an ordered partition $(A_1\LL A_{\ell})$ an {\em $\ell$-ring}. (In fact \ref{lring} below says a strongly-connected digraph admits 
an $\ell$-ring if and only if every directed cycle has length a multiple of $\ell$.) 

\begin{figure}[ht]
\centering

\begin{tikzpicture}[auto=left]

\draw[rotate=60] (2,0) ellipse (10pt and 6pt);
\draw[rotate=30] (2,0) ellipse (10pt and 6pt);
\draw[] (2,0) ellipse (10pt and 6pt);
\draw[rotate=-30] (2,0) ellipse (10pt and 6pt);
\draw[rotate=-60] (2,0) ellipse (10pt and 6pt);

\draw[domain=70:290,smooth,variable=\x, thick, dotted] plot ({2*cos(\x)},{2*sin(\x)});

\def\s{1.8}
\def\t{2.2}
\def\r{.98}

\draw[-] ({\s*cos(57)},{\s*sin(57)})-- ({\s*cos (33)},{\s*sin(33)});
\draw[-] ({\t*cos(57)},{\t*sin(57)})-- ({\t*cos (33)},{\t*sin(33)});
\draw [-{Stealth[scale=1.2]}] ({\s*cos (33)},{\s*sin(33)})--({\r*\s*cos (50)},{\r*\s*sin(50)});
\draw [-{Stealth[scale=1.2]}] ({\t*cos (33)},{\t*sin(33)})--({\r*\t*cos (50)},{\r*\t*sin(50)});

\draw[-] ({\s*cos(27)},{\s*sin(27)})-- ({\s*cos (3)},{\s*sin(3)});
\draw[-] ({\t*cos(27)},{\t*sin(27)})-- ({\t*cos (3)},{\t*sin(3)});
\draw [-{Stealth[scale=1.2]}] ({\s*cos (3)},{\s*sin(3)})--({\r*\s*cos (20)},{\r*\s*sin(20)});
\draw [-{Stealth[scale=1.2]}] ({\t*cos (3)},{\t*sin(3)})--({\r*\t*cos (20)},{\r*\t*sin(20)});

\draw[-] ({\s*cos(-3)},{\s*sin(-3)})-- ({\s*cos (-27)},{\s*sin(-27)});
\draw[-] ({\t*cos(-3)},{\t*sin(-3)})-- ({\t*cos (-27)},{\t*sin(-27)});
\draw [-{Stealth[scale=1.2]}] ({\s*cos (-27)},{\s*sin(-27)})--({\r*\s*cos (-10)},{\r*\s*sin(-10)});
\draw [-{Stealth[scale=1.2]}] ({\t*cos (-27)},{\t*sin(-27)})--({\r*\t*cos (-10)},{\r*\t*sin(-10)});

\draw[-] ({\s*cos(-33)},{\s*sin(-33)})-- ({\s*cos (-57)},{\s*sin(-57)});
\draw[-] ({\t*cos(-33)},{\t*sin(-33)})-- ({\t*cos (-57)},{\t*sin(-57)});
\draw [-{Stealth[scale=1.2]}] ({\s*cos (-57)},{\s*sin(-57)})--({\r*\s*cos (-40)},{\r*\s*sin(-40)});
\draw [-{Stealth[scale=1.2]}] ({\t*cos (-57)},{\t*sin(-57)})--({\r*\t*cos (-40)},{\r*\t*sin(-40)});

\end{tikzpicture}

\caption{An $\ell$-ring.} \label{fig:ring}

\end{figure}
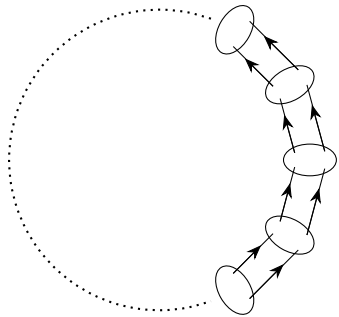
So we can confine ourselves to digraphs that 
admit an $\ell$-ring, and consequently all directed cycles have length divisible by $\ell$. But admitting an $\ell$-ring is not enough 
to guarantee that $G$ is $\ell$-cyclic; there might still be 
directed cycles with winding number more than one, in the natural sense, and we 
need to put additional restrictions on the $\ell$-ring. There are two restrictions we want to consider.

First, let us say that
a digraph $G$, with an $\ell$-ring $(A_1\LL A_{\ell})$, admits an  {\em $\ell$-annular drawing} if $G$ can be drawn without crossings 
in the plane $\mathbb{R}^2$, such that for $1\le i\le \ell$,
the vertices in $A_i$ are represented by points $(x,y)$ in the half-line 
$$\{(r \cos(2\pi i/\ell), r \sin(2\pi i/\ell)): r>0\},$$ 
and the edges from $A_{i-1}$ to $A_{i}$ are drawn within the area
$$\{(r\cos \theta, r\sin \theta): r>0, 2\pi (i-1)/\ell\le \theta\le 2\pi i/\ell\}.$$
We say $G$ is {\em $\ell$-annular} if $G$ admits an $\ell$-annular drawing.
It is easy to see that $\ell$-annular digraphs are $\ell$-cyclic, because closed curves in the plane that wind more than once around the origin
must intersect themselves.
\begin{figure}[ht]
\centering

\begin{tikzpicture}[auto=left]
\tikzstyle{every node}=[inner sep=1.5pt, fill=black,circle,draw]

\node (a1) at (1,0) {};
\node (a2) at (2,0) {};
\node (a3) at (3,0) {};
\node (b1) at (0,1) {};
\node (b2) at (0,2) {};
\node (c1) at (-1,0) {};
\node (c2) at (-2,0) {};
\node (c3) at (-3,0) {};
\node (d1) at (0,-1) {};
\node (d2) at (0,-2) {};
\node (d3) at (0,-3) {};

\foreach \from/\to in {a1/b1,a2/b1,a2/b2,a3/b2, b1/c1, b1/c2, b2/c2, b2/c3, c1/d1,c1/d2,c2/d2,c2/d3,c3/d3,d1/a1,d2/a1,d2/a2,d2/a3,d3/a3}
\draw [] (\from)--(\to);

\def\s{3/2}
\tikzstyle{every node}=[];
\draw [-{Stealth[scale=\s]}] (a1) -- (.5,.5);
\draw [-{Stealth[scale=\s]}] (a2) -- (1,.5);
\draw [-{Stealth[scale=\s]}] (a2) -- (1,1);
\draw [-{Stealth[scale=\s]}] (a3) -- (1.5,1);
\draw [-{Stealth[scale=\s]}] (b1) -- (-.5,.5);
\draw [-{Stealth[scale=\s]}] (b1) -- (-1,.5);
\draw [-{Stealth[scale=\s]}] (b2) -- (-1,1);
\draw [-{Stealth[scale=\s]}] (b2) -- (-1.5,1);
\draw [-{Stealth[scale=\s]}] (c1) -- (-.5,-.5);
\draw [-{Stealth[scale=\s]}] (c1) -- (-.5,-1);
\draw [-{Stealth[scale=\s]}] (c2) -- (-1,-1);
\draw [-{Stealth[scale=\s]}] (c2) -- (-1,-1.5);
\draw [-{Stealth[scale=\s]}] (c3) -- (-1.5,-1.5);
\draw [-{Stealth[scale=\s]}] (d1) -- (.5,-.5);
\draw [-{Stealth[scale=\s]}] (d2) -- (.5,-1);
\draw [-{Stealth[scale=\s]}] (d2) -- (1,-1);
\draw [-{Stealth[scale=\s]}] (d2) -- (1.5,-1);
\draw [-{Stealth[scale=\s]}] (d3) -- (1.5,-1.5);

\end{tikzpicture}

\caption{A $4$-annular drawing.} \label{fig:annular}

\end{figure}
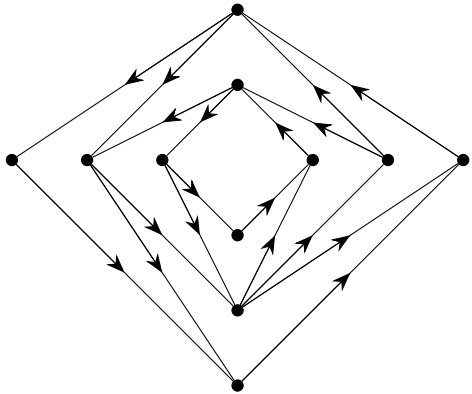

Second, a more brutal restriction:
let us say that an $\ell$-ring $(A_1\LL A_{\ell})$ is {\em pinched} if some $A_i$ has cardinality one, and $G$ is {\em $\ell$-pinched} if it admits a pinched $\ell$-ring.
Every $\ell$-pinched digraph is $\ell$-cyclic. 

One might hope for a construction for all $\ell$-cyclic digraphs, somehow building them out of classes of $\ell$-cyclic digraphs 
that we understand; and the two above seem to be natural classes to use as building blocks. There are theorems that say that some sort 
of planar drawing is equivalent to the nonexistence of certain pairs of disjoint paths, in graphs~\cite{GM9, 2paths, thomassenpaths} 
and in acyclic digraphs~\cite{acyclic,thomassenacyclic}, so at first sight one might hope for something similar.
We do not know such a construction when $\ell\ne 3$, but when $\ell=3$ we can do it. (Rather sadly, the first, more interesting, 
class above is not involved.) We need a
composition operation. Let us say an edge $uv$ of a digraph $G$ is {\em special} if either $u$ has outdegree one, or
$v$ has indegree one. If $G,G'$ are $\ell$-cyclic, and $uv$ is a special edge of $G$, and $u'v'$ is a special edge of $G'$, we can make the
idenifications $u=u'$ and $v=v'$, and the digraph we produce is also $\ell$-cyclic. Let us call this operation a {\em special edge sum}.
We will prove:
\begin{thm}\label{all3cycled}
Every unbreakable 3-cyclic digraph can be made by special edge sums, starting from 3-pinched digraphs.
\end{thm}

If $G$ is a graph or digraph and $X\subseteq VG)$, $G[X]$ denotes the subgraph or subdigraph induced on $X$.
Let us say a {\em tricycle} is a directed cycle of length three.
Take a graph consisting of a cycle $C$                             
of even length, and an extra vertex $v$ adjacent to every vertex of $C$; and now direct its edges, making a digraph $G$, so that 
every edge of $C$ makes a tricycle with $v$.  (Consequently the edges of $C$ are directed alternately clockwise and counterclockwise.)
We call the digraph we obtain a {\em diwheel}, and $v$ is its {\em hub}. It is easy to see that diwheels are 3-cyclic, and are not 3-annular.                  
We say $G$ is {\em diwheel-free} if no subdigraph is a diwheel. 
We will see that general unbreakable 3-cyclic digraphs can be made by piecing together on small cutsets 
digraphs that are either 3-pinched or diwheel-free, so let us try to understand those that are diwheel-free.

If $G$ is 3-cyclic, there are several equivalent descriptions of when $G$ is diwheel-free: 
\begin{thm}\label{nowheel}
Let $G$ be an unbreakable 3-cyclic digraph. Then the following are equivalent:
\begin{itemize}
\item $G$ is diwheel-free;
\item $G$ admits a 3-ring $(A_1,A_2,A_3)$ such that each of the graphs $G^-[A_1\cup A_2]$, $G^-[A_2\cup A_3]$, $G^-[A_3\cup A_1]$ is a tree;
\item $G$ can be built starting from a tricycle by repeatedly adding a new vertex with exactly one out-neighbour and exactly one in-neighbour;
\item $G$ has exactly $2|G|-3$ edges.
\end{itemize}
\end{thm}
``3-cyclic'' is part of the hypothesis in the above; being buildable in the way described does not guarantee that $G$ is 3-cyclic.
But here is one that does. Let us say that a digraph $G$ is {\em safely buildable} if $G$ can be constructed by starting with a tricycle,
and repeating the following
operation:
given a digraph $H$, choose a vertex $v$ of $H$ with degree 2 (so in-degree and out-degree both one);
choose one of its neighbours $u$; and add a set of new vertices, each of degree two and making a tricycle with $u,v$. 

\begin{thm}\label{safeadd}
For a digraph $G$, the following are equivalent:
\begin{itemize}
\item $G$ is unbreakable, 3-cyclic and diwheel-free;
\item $G$ is safely buildable.
\end{itemize}
\end{thm}

It is tempting to conjecture that all 
unbreakable 3-cyclic diwheel-free digraphs are 3-annular, but that is false. The digraphs in Figure 
\ref{fig:branchers} are unbreakable and 3-cyclic, but not 3-annular. We call these two digraphs, together with 
the two further digraphs obtained by reversing all edges of these two, the four {\em branchers}.
We will prove:
\begin{thm}\label{planar}
If $G$ is an unbreakable digraph, then $G$ is 3-annular if and only if it contains no diwheel and none of the four branchers.
\end{thm}

We conclude this section by proving the statement above about $\ell$-rings:
\begin{thm}\label{lring}
Let $\ell\ge 2$ be an integer. If $G$ is a strongly connected digraph, then every directed cycle of $G$ has length a multiple of $\ell$ if and only if $G$ admits an $\ell$-ring.
\end{thm}
\Proof
The ``if'' part is clear, and we must prove ``only if''. Suppose then that $G$ is strongly connected and every directed cycle has length a multiple of $\ell$.
It follows that the length of every directed closed walk is also a multiple of $\ell$ (for instance, by induction on the length of 
the walk: if it has a ``repeated'' vertex, its length is the sum of the lengths of two shorted directed closed walks, and otherwise 
it forms a directed cycle). Fix a vertex $r$. We claim that if $W_1,W_2$ are directed walks from $r$ to a vertex $v$, then their lengths 
are equal modulo $\ell$; because there is a directed walk  $W$ from $v$ to $r$, and the directed closed walks formed by concatenating $W_1, W$ 
and concatenating $W_2,W$ both have length zero modulo $\ell$. For $0\le i\le \ell-1$, let $A_i$ be the set of vertices $v$ such that all
directed walks from $r$ to $v$ have length $i$ modulo $\ell$; then $(A_0\LL A_{\ell-1})$ is an $\ell$-ring. This proves \ref{lring}.~\bbox

\section{Connectivity in a 3-cyclic digraph}

If $G$ is $\ell$-cyclic and strongly connected, then for every edge $uv$, since there is a directed path from $v$ to $u$ and its union 
with $uv$ is a directed cycle, it follows that $uv$ belongs to a directed cycle of length $\ell$. In particular, every edge of a 
strongly connected  3-cyclic digraph belongs to a tricycle. (We will use this fact without further explanation.)
We have nothing else to say about $\ell$-cyclic digraphs in general, and from now on we just consider the case $\ell=3$. 
In this section we prove some lemmas about connectivity in $G^-$, when $G$ is 3-cyclic.
If $P$ is a path (not necessarily a directed path) of a digraph $G$, with ends $a,b$ say, we say an edge $e=uv$ of $P$ is 
{\em directed towards $b$ in $P$} if $v,b$ belong to the same component of $P^-\setminus e$, and otherwise it is {\em directed away from $b$}.
We will prove the following useful lemma:
\begin{thm}\label{cutset}
Let $H$ be a strongly connected 
subdigraph of a strongly connected 3-cyclic digraph $G$. Suppose that $X\subseteq V(H)$, 
and $(A,B)$ is a partition of $V(H)\setminus X$ such that $A,B\ne \emptyset$, and there are no edges of $H$ from $B$ to $A$.
Suppose that
$P$ is a path of $G\setminus X$ (not necessarily directed) between $A,B$, such that as few edges as possible are directed away from $b$ in $P$, and subject to that, the length of $P$ is minimum.
Then either:
\begin{itemize}
\item there is a tricycle with vertices $u,v,w$, such that $u,v\in X$ and $w$ belongs to the interior of $P$; or
\item there exists $x\in X$ such that every vertex of $P$ is adjacent to or from $x$, and every edge of $P$ forms a tricycle with $x$.
\end{itemize}
\end{thm}
\Proof From the minimality of the length of $P$, it follows that only the first vertex of $P$ (say $a$) is in $A$, and only its last vertex
(say $b$) is in $B$. Suppose first that no edges of $P$ are directed away from $b$ in $P$, and so $P$ is a directed path from $a$ to $b$.
Since $H$ is strongly-connected, there is a directed path $Q$ of $H$ from $b$ to $a$, with length at least two, since
there are no edges of $H$ from $B$ to $A$. But $P\cup Q$ is a directed cycle, and hence has length three; so $P$ has length one, and 
$Q$ has length two, with middle vertex in $X$, and the theorem holds.

\bigskip

Let $F$ be the set of edges of $P$ directed away from $b$; so we may assume that $F\ne \emptyset$. For
each $uv\in X$, we say $x\in X$ is an {\em apex} for $uv$ if 
$u\rightarrow v\rightarrow x\rightarrow u$
is a directed cycle.
\\
\\
(1) {\em Each edge in $F$ has a unique apex.}
\\
\\
Let $uv\in F$. Since $G$ is strongly connected,
there is a directed path in $G$ from $v$ to $u$, which makes a directed cycle with the edge $uv$, and so has length two. Let $v\rightarrow  w\rightarrow u$
be this directed path. It follows that $w\notin V(P)$, from the minimality of the length of $P$. If $w\in A$, then 
the path $w\rightarrow u$ togther with the subpath of $P$ between $u,b$ makes a path of $G\setminus X$ between $A,B$, and it has 
fewer edges directed away from $b$, a contradiction. Similarly $w\notin B$. If $w\notin X$, then $w\notin V(H)$, and the union of 
$v\rightarrow w \rightarrow u$ and the subpaths of $P$ between $a,v$ and between $u,b$ is a path of $G\setminus X$ between $A,B$, and it has  
fewer edges directed away from $b$, a contradiction. So $w\in X$, and hence $w$ is an apex for $uv$. Suppose that $w'\in X$ is another
apex for $uv$ with $w'\ne w$. Choose a directed path $Q$ of $H$ from $w$ to $w'$. Then the union of $Q$ with the path 
$w'\rightarrow u\rightarrow v\rightarrow$ makes a directed cycle of length more than three, a contradiction. This proves (1).
\\
\\
(2) {\em If $F$ contains two edges with a common end, say $uv,vw\in F$, then the first bullet of the theorem holds.}
\\
\\
Let $x,x'$ be apexes for $uv$ and for $vw$ respectively. Choose a directed path $Q$ of $H$ from $x$ to $x'$. Then the union of $Q$
with $x'\rightarrow v\rightarrow x$ is a directed cycle, and so $Q$ has length one, and the first bullet holds. This proves (2).

\bigskip

In view of (2) we may assume that no two edges in $F$ have a common end. Let the vertices of $P$ in order be
$a=v_1\DD v_2\CC v_k=b$, and let $F=\{v_{i+1}v_i:(i\in I)\}$. Thus $I\subseteq \{1\LL k-1\}$, and every two members of $I$ differ 
by at least two.
\\
\\
(3) {\em If $i,j\in I$ such that $j>i$ and none of $i+1\LL j-1$ belong to $F$, then $j=i+2$
and $v_{i+1}v_i$ and $v_{j+1}v_j$ have the same apex.}
\\
\\
Let $v_{i+1}v_i$ and $v_{j+1}v_j$
have apexes $x,x'$ respectively, and let $Q$ be a directed path of $H$ from $x'$ to $x$. Let $R$ be the subpath of $P$ from $v_{i+1}$
t $v_j$; thus $R$ is a directed path from $v_{i+1}$ to $v_j$, with length at least one. But the union of $Q,R$ and the edges $xv_{i+1}$ and $v_jx'$
makes a directed cycle, and so $R$ has length one and $Q$ has length zero. This proves (3). 

\bigskip

From (3) it follows that all edges in $F$ have the same apex $x\in X$. Let $h,j\in F$ be the smallest and largest members of $F$.
From (3) it also follows that every vertex of the subpath of $P$ between $v_h$ and $v_{j+1}$ is adjacent to or from $x$, and each edge of
this subpath makes a directed cycle with $x$. We claim the same is true for the initial subpath $v_1\DD v_h$ and terminal 
subpath $v_{j+1}\DD v_k$. To see that every edge of the subpath $v_1\DD v_h$ makes a tricycle with $x$, we may assume that $h>1$.
Let $Q$ be a directed path of $H$ from $x$ to $a$; then the union of $Q$, the subpath of $P$ from $v_1=a$ to $v_h$,
and the edge $v_hx$, makes a directed cycle, which must have length three. So $h=2$ and $v_1,v_2$ are adjacent from $x$ and to $x$
respectively. Similarly, either $j+1=k$, or $j+1=k-1$ and $v_{k-1}, v_k$ are adjacent from $x$ and to $x$ respectively. 
Consequently the second bullet of the theorem holds. This proves \ref{cutset}.~\bbox

Here is another lemma we will need later, on ``ear decompositions''. The following ``ear'' results are well known:
\begin{itemize}
\item If $H$ is a non-null subdigraph of a strongly connected digraph $G$, and $H\ne G$, either there is a directed path of $G$
of length at least one,  with both ends in $V(H)$ and no edge or internal vertex in $H$, or there is a directed cycle with exactly one vertex in $H$
\item If $H$ is a subgraph with $|H|\ge 2$ of a 2-connected graph $G$, and $H\ne G$, there is a path of $G$ of length at least one, 
with both ends in $V(H)$ and no edge or internal vertex in $V(H)$.
\end{itemize}
We want to make a sort of common generalization. If $H$ is a subdigraph of a digraph $G$, let us say an {\em ear} (for $H$, in $G$) is a directed path of 
$G$ of length at least one, with both ends in $V(H)$ and with no edge or internal vertex in $H$.
\begin{thm}\label{ears}
Let $H$ be a subdigraph with $|H|\ge 2$ of an unbreakable digraph $G$, with $H\ne G$.
Then there is an ear.
\end{thm}
\Proof
Let $G'$ be obtained from $G$ by deleting all edges of $H$. Each vertex $v$ of $G$ belongs to a unique strong component 
$K_v$ of $G'$. If $K_u=K_v$ for some distinct $u,v\in V(H)$, then a directed path of $K_v$ from $u$ to $v$ is an ear, so we assume that
$K_v\;(v\in V(H))$ are all different. Suppose that there is a strong component $K$ of $G'$ that contains no vertex of $H$.
Choose $v\in V(K)$, and let $P,Q$ be minimal directed paths of $G$ from $v$ to $V(H)$ and from $V(H)$ to $v$ respectively. Consequently $P,Q$ are both paths of $G'$.
Let $P$ be from $v$ to $w\in V(H)$, and let $Q$ be from $u\in V(H)$ to $v$. If $u=w$, then $P\cup Q$ is strongly connected, and so 
is a subdigraph of $K$; but then $K$ contains $u$, a contradiction. If $u\ne w$, there is a directed path in $P\cup Q$ from $u$ to $w$,
and this is an ear. So there is no such strong component $K$, and hence every strong component of $G'$ contains exactly one vertex of $H$.

For each $v\in V(H)$, we may assume that there is no edge $xy$ between $V(K_v)\setminus \{v\}$ and $V(H)\setminus \{v\}$; because if $x\in V(K_v)$, then 
the union of a directed path of $K_v$ from $v$ to $x$ and the edge $xy$ is an ear, while if $y\in V(K_v)$ then the union of the edge $xy$
and a directed path of $K_v$ from $y$ to $v$ is an ear.
Similarly, if $u,v\in V(H)$ are distinct, we may assume that there is no edge $xy$ from $V(K_u)\setminus \{u\}$ to $V(K_v)\setminus \{v\}$,
since the union of a directed path of $K_u$ from $u$ to $x$, the edge $xy$, and a directed path of $K_v$ from $y$ to $v$ would be an ear.

It follows that if $v\in V(H)$, there is no path of $G^-$ between $V(K_v)\setminus \{v\}$ and $V(G)\setminus V(K_v)$.
Since $|H|\ge 2$ and $G$ is weakly 2-connected, it follows that $K_v=\{v\}$ for each $v\in V(H)$. Since every vertex of $G$ belongs
to $G'$ and so to some $K_v$, it follows that $V(H)=V(G)$. Consequently any edge of $G$ not in $E(H)$ is an ear, and so we may
assume there is no such edge; and so $G=H$, a contradiction. This proves \ref{ears}.~\bbox

The point of the theorem is, if $H$ is an unbreakable subdigraph of an unbreakable digraph $G$, and we add an ear to $H$, we obtain another unbreakable subdigraph, and we can repeat the process 
until we reach $G$.

\section{3-cyclic digraphs with diwheels}
First, let us prove the following, a step towards \ref{all3cycled}. 
\begin{thm}\label{block2}
If $G$ is an unbreakable 3-cyclic digraph, containing a diwheel with hub $v$, then either:
\begin{itemize}
\item $G$ admits a $3$-ring $(A_1,A_2,A_3)$ with $|A_1|=\{v\}$; or
\item there is a diwheel $W$ with hub $v$, and a vertex $u\in V(W)\setminus \{v\}$,  such that $\{u,v\}$ is a cutset of $G^-$.
\end{itemize}
\end{thm}
\Proof
Choose a set $A_2$ of out-neighbours of $v$, and a set $A_3$ of in-neighbours of $v$, such that $G[A_2\cup A_3]$ is weakly 2-connected,
with $A_2\cup A_3$ maximal. (This is possible since there is a diwheel with hub $v$.) By \ref{lring}, $G$ admits a 3-ring $(B_1,B_2,B_3)$ with $v\in B_1$; and so $A_2\subseteq B_2$ and $A_3\subseteq B_3$.
If $A_2\cup A_3=V(G)\setminus \{v\}$, then 
$B_1=\{v\}$, and the first bullet holds; so we may assume that $A_2\cup A_3\ne V(G)\setminus \{v\}$. 
Since $G[A_2\cup A_3]$ is weakly 2-connected and every edge with both ends in $B_2\cup B_3$ is from $B_2$ to $B_3$, it follows that each vertex in $A_2$ has
an out-neighbour in $A_3$ and vice versa; and so every vertex in $A_2\cup A_3$ is contained in a cyclic triangle that also contains $v$. 
Consequently $G[A_2\cup A_3\cup \{v\}]$ is strongly connected. By \ref{ears}, there is an ear $Q$ for $G[A_2\cup A_3\cup \{v\}]$
in $G$. Let $Q$ be from $a$ to $b$. Then $Q$ has length at least two, since $G[A_2\cup A_3\cup \{v\}]$ is an induced subdigraph. 
Since there is a directed path of $G[A_2\cup A_3\cup \{v\}]$ from $b$ to $a$, and its union with $Q$ is a directed cycle, it follows that
$Q$ has length two, and $a$ is adjacent from $b$. Let $c$ be the middle vertex of $Q$. Suppose that $a,b\ne v$; then $a,b\in A_2\cup A_3$,
and so $a\in A_3$ 
and $b\in A_2$, since $ba$ is an edge. Since $G[A\cup B]$ is weakly 2-connected, there exists $b'\ne b$ in $A_2$ adjacent to $a$, and 
there exists $a'\ne a$ in $A_3$ adjacent from $b$. But then 
$$v\rightarrow b'\rightarrow a\rightarrow c\rightarrow b\rightarrow a'\rightarrow v$$ 
is a directed cycle, a contradiction.

So one of $a,b$ equals $v$, and from the symmetry under reversing the direction of all edges, we may assume that $b=v$, and consequently
$a\in A_2$. Since $G[A_2\cup A_3]$ is weakly 2-connected, it has a cycle containing $a$, and so there is a diwheel with hub $v$,
containing $a$, and with vertex set a subset of $A_2\cup A_3\cup \{v\}$. Thus if $\{v,a\}$ is a cutset of $G^-$,
the second bullet of the theorem holds. Let $C$ be the set of all vertices of $G$ that are not in $A_3$, and are adjacent from $a$ and to $v$
(thus $c\in C$), and let $A=A_2\cup A_3\setminus \{a\}$. Since we may assume that $\{v,a\}$ is not a cutset of $G^-$, 
there is a path $P$ of $G^-$ between $A$ and $C$; and so we may choose $P$ as in \ref{cutset}. No vertex of the interior of $P$ belongs to
$A\cup C\cup \{a,v\}$, and so no vertex of the interior of $P$ is adjacent to $v$ and from $a$, and so the first bullet of \ref{cutset}
does not hold. Consequently the second holds, and there exists $x\in \{v,a\}$ adjacent to or from every vertex of $P$, such that every edge of $P$ makes a tricycle with $x$. Since $G[A_2\cup A_3\cup V(P)]$ is weakly 2-connected, the maximality of $A_2\cup A_3$
implies that some vertex of $P$ is not a neighbour of $v$, and so $x\ne v$. Hence $x=a$, and so the vertex $p$ of $P$ in $A_2\cup A_3$
belongs to $A_3$. Let $q$ be the neighbour of $p$ in $P$ (so $p\rightarrow q \rightarrow a$ is a directed path). Choose $a'\in A_2\setminus \{a\}$ adjacent to $p$; and choose $p'\in A_3\setminus \{p\}$ adjacent from $a$. Then 
$$v\rightarrow a'\rightarrow p\rightarrow q\rightarrow a\rightarrow p'\rightarrow v$$
is a directed cycle, a contradiction. This proves \ref{block2}.~\bbox

Now let us convert this into a construction for all 3-cyclic digraphs.
We recall that an edge $uv$ of a digraph $G$ is {\em special} if either $u$ has outdegree one, or
$v$ has indegree one. 

\begin{thm}\label{specialcomp}
Let $G$ be an unbreakable, 3-cyclic digraph, and let $\{u,v\}$ be a cutset of $G^-$.
Then $u,v$ are adjacent in $G$, and we assume that $uv\in E(G)$ without loss of generality. For every component $B$ of $G^-\setminus \{u,v\}$, $uv$ is a special edge of $G[B\cup \{u,v\}]$.
\end{thm}
\Proof
Since $\{u,v\}$ is a cutset of $G^-$, there is a partition $(L,R)$ of $V(G)\setminus \{u,v\}$ with $L,R\ne \emptyset$, 
such that there are no edges of $G^-$ between $L,R$. By exchanging $u,v$ if necessary, we may assume that $vu\notin E(G)$. 
Since $G$ is strongly connected, there is a directed path $P$ from $v$ to $u$. 
We may assume that $P$ is a path of $G[L\cup \{u,v\}]$. 
\\
\\
(1) {\em $uv\in E(G)$, and there is a tricycle containing $uv$ with third vertex in $L$.}
\\
\\
We claim first that there is a directed path of
 $G[L\cup \{u,v\}]$ from $u$ to $v$. This is true if $uv\in E(G)$; and otherwise, since there is a directed path from $u$ to $v$ in $G$,
and its union with $P$ cannot be a directed cycle, this path is also in  $G[L\cup \{u,v\}]$. This proves that 
there is a directed path $Q$ of
 $G[L\cup \{u,v\}]$ from $u$ to $v$. Hence $P\cup Q$ is strongly connected, and so $G[L\cup \{u,v\}]$ is strongly connected.
There is a path of $G^-$ between $u,v$ disjoint from $L$. Let us apply \ref{cutset}, taking $H=G[L\cup \{u,v\}]$, $X=L$, $A=\{u\}$
and $B=\{v\}$. The first outcome is impossible, since there are no edges between $L,R$, and so the second holds, that is, 
there exists $x\in L$ and a path $R$ of $G^-\setminus L$ between $u,v$, such that every vertex of $R$ is adjacent to or from $x$, 
and every edge of $P$ forms a tricycle with $x$. Consequently $R$ has no internal vertex, and so has length one; and so $uv\in E(G)$ (since $vu\notin E(G)$), and the edge $uv$ forms a triangle with $x$. This proves (1).
\\
\\
(2) {\em There is a tricycle containing $uv$ with third vertex in $R$.}
\\
\\
There is no directed path of $G[R\cup \{u,v\}]$ from $u$ to $v$ not using the edge $uv$, since its union with $P$
would be a directed cycle of length at least four. 
Choose a path $S$ of $G[R\cup \{u,v\}]$ between $u,v$ not using the edge $uv$, and with as few edges directed away from $v$ as possible.
It has at least one edge directed away from $v$, say
$xy$. Let $xyz$ be a tricycle containing $xy$. Rerouting $S$ from $y$ to $x$ along $y\rightarrow z\rightarrow x$
does not give a path of $G[R\cup \{u,v\}$ between $u,v$ with fewer edges directed away from $v$ and not using the edge $uv$; so
the edge $uv$ belongs to the tricycle, and so there is a tricycle of $G[R\cup \{u,v\}$ containing $uv$.
This proves (2).

\bigskip

From (1) and (2), every directed path of $G$ from $u$ to $v$ uses the edge $uv$. Consequently there is a partition $(X,Y)$
of $V(G)$ with $u\in X$ and $v\in Y$, such that $uv$ is the only edge of $G$ from $X$ to $Y$. If $yx\in E(G)$,
and $x\in X$ and $y\in Y$, then since $xy$ belongs to a tricycle, that tricycle must contain $uv$, and so either $x=u$ or $y=v$.
Now $B$ (in the theorem) is a component of $G^-\setminus \{u,v\}$, and yet there is no edge of $B$ between $B\cap X$ and $B\cap Y$; and so $B$
is a subset of one of $X,Y$. But then $uv$ is a special edge of $G[B\cup \{u,v\}]$. This proves \ref{specialcomp}.~\bbox

If $G,G'$ are digraphs, and $uv$ is a special edge of $G$, and $u'v'$ is a special edge of $G'$, and $G''$ is obtained 
by making the identifications $u=u'$ and $v=v'$, we say that $G''$ is a {\em special edge sum} of $G,G'$. 
\begin{thm}\label{findspecial}
Let $G$ be an unbreakable 3-cyclic digraph, and let $uv\in E(G)$ such that $\{u,v\}$ is a cutset of $G^-$. 
There is a partition of $V(G)\setminus \{u,v\}$ into $A,B$ say, such that 
\begin{itemize}
\item $A,B\ne \emptyset$; 
\item there are no edges between $A,B$;
\item $uv$ is a special edge of each of $G[A\cup \{u,v\}], G[B\cup \{u,v\}]$; and
\item $G[A\cup \{u,v\}], G[B\cup \{u,v\}]$ are unbreakable.
\end{itemize}
\end{thm}
\Proof
We claim first that there exist $A,B$ satisfying the first three bullets.
Let $A_1\LL A_m,B_1\LL B_n$ be the components of $G^-\setminus \{u,v\}$, where $A_1\LL A_m$ are the components that contain 
an out-neighbour of $u$. Thus $m+n\ge 2$.  
By \ref{specialcomp}, $A_1\LL A_m$ contains no in-neighbours of $v$, and so $uv$ is a special edge
of $G[A_1\cupcup A_m\cup \{u,v\}]$. Moreover, $uv$ is a special edge of $G[B_1\cupcup B_n\cup\{u,v\}]$ since $B_1\cupcup B_n$
contains no out-neighbour of $v$; and so if $m,n>0$ then the claim is satisfied by setting $A=A_1\cupcup A_m$ and $B=B_1\cupcup B_n$.
Hence we assume that $m=0$ without loss of generality. But then $uv$ is a special edge of $G$, and the claim is satisfied by setting $A=B_1$ and $B=B_2\cupcup B_n$. This proves that there exist $A,B$ satisfying the first three bullets.

We claim that the fourth bullet of the theorem is also satisfied. Clearly $G[A\cup \{u,v\}], G[B\cup \{u,v\}]$
are weakly 2-connected, and we just need to show they are strongly connected. From the symmetry, it suffices to show that 
$G[A\cup \{u,v\}]$ is strongly connected, and hence we only need to show there is a directed path in $G[A\cup \{u,v\}]$ from $v$ to $u$,
since $G$ is strongly connected. Again, without loss of generality we may assume that $v$ is the only out-neighbour of $u$ in 
$G[A\cup \{u,v\}]$. Since $A\ne \emptyset$ and $G$ is weakly 2-connected, $u$ has a neighbour in $A$, say $a$; and hence $au\in E(G)$.
There is a directed path in $G$ from $v$ to $a$, and it uses no edge between $u,A$; so it is a path
of $G[A\cup \{u,v\}]$. Adding the edge $au$ gives the directed path we need. This proves \ref{findspecial}.~\bbox

That takes us a step closer to \ref{all3cycled}:
\begin{thm}\label{all3cycled2}
Every unbreakable 3-cyclic digraph can be made by special edge sums, starting from 3-pinched digraphs
and unbreakable 3-cyclic diwheel-free digraphs.
\end{thm}
\Proof Let $G$ be an unbreakable 3-cyclic digraph; we prove by induction on $|G|$ that $G$
can be made by special edge sums, starting from 3-pinched digraphs and unbreakable 3-cyclic diwheel-free digraphs. If $G$ is 3-pinched
or diwheel-free, this is true, and otherwise, by 
\ref{block2}, $G^-$ is not 3-connected. From \ref{findspecial}, $G$ is expressible as a special edge sum of two smaller 
unbreakable 3-cyclic digraphs, and the result follows from the inductive hypothesis.~\bbox

\section{3-cyclic digraphs without diwheels}
In this section we prove \ref{safeadd}, and complete the proof of \ref{all3cycled}, but we need a few lemmas. 

\begin{thm}\label{nbrs}
Let $G$ be an unbreakable 3-cyclic digraph, let $v\in V(G)$, and let $N$ be the set of vertices that are
adjacent with $v$ in $G^-$. Then $G[N]$ is weakly connected.
\end{thm}
\Proof
Suppose that there is a partition of $N$ into two nonempty sets $A,B$, and there is no edge of $G^-$ between $A,B$.
For each $u\in A$, the edge between $u,v$ is contained in a tricycle (because $G$ is strongly connected and 3-cyclic), and all
three of its vertices belong to $A\cup \{v\}$ (because there is not edge of $G^-$ between $A,B$). Consequently
$G[A\cup \{v\}]$ is strongly connected, and similarly so is $G[B\cup \{v\}]$, and hence so is $G[N\cup \{v\}]$.
There is a path of $G^-$ between $A,B$, since $G$ is weakly 2-connected; and so by
\ref{cutset}, applied with $H=G[N\cap \{v\}]$ and $X=\{v\}$, there is such a path such that all its vertices belong to $N$,
a contradiction. This proves \ref{nbrs}.~\bbox

\begin{thm}\label{diamond}
Let $G$ be an unbreakable, 3-cyclic, diwheel-free digraph. Let $xy\in E(G)$, and let both $u,v\in V(G)$
make tricycles with the edge $xy$. Then there is no path of $G^-\setminus \{x,y\}$ between $u,v$.
\end{thm}
\Proof Let $A$ be the set of all vertices that are adjacent to $x$ and from $y$ (so $u,v\in A$). We will prove that all vertices in $A$
belong to different components of $G^-\setminus \{x,y\}$. Suppose not, and choose $b\in A$ such that there is a path of $G^-\setminus \{x,y\}$
between $b$ and some other vertex in $A$ . Let $P$ be a path of $G^-\setminus \{x,y\}$, with ends $a,b\in A$ say, 
with $k$ edges of $P$ directed away from $b$,
chosen such that:
\begin{itemize}
\item every path of $G^-\setminus \{x,y\}$ between $b$ and $A\setminus \{b\}$ has at least $k$ edges directed away from $b$; and
\item among all paths of $G^-\setminus \{x,y\}$ between $a,b$ that have only $k$ edges directed away from $b$, $P$ has minimum length.
\end{itemize}
Now $ab\notin E(G)$ since $a\rightarrow b\rightarrow x\rightarrow y\rightarrow a$
is not a directed cycle, and similarly $ba\notin E(G)$. Moreover, from the second bullet above, no vertex of the interior of $P$
makes a tricycle with $xy$ (because such a vertex would be in $A$, and we could choose $P$ shorter). 
By \ref{cutset} applied with $X=\{x,y\}$ and $H=G[\{a,b,x,y\}]$, there exists $z\in \{x,y\}$ such that every edge of $P$
makes a tricycle with $z$. But then $G$ contains a diwheel with hub $z$, a contradiction. This proves \ref{diamond}.~\bbox

\begin{thm}\label{oneper}
Let $G$ be an unbreakable, 3-cyclic, diwheel-free digraph, let $T$ be a tricycle, and let $B$ be a
component  
of $G^-\setminus V(T)$. Exactly two of the vertices of $T$ have a neighbour in $B$, say $t_1,t_2$; and there is a unique $b\in B$ such that
$\{t_1,t_2,b\}$ is the vertex set of a tricycle.
\end{thm}
\Proof Let $T$ be $t_1\rightarrow t_2\rightarrow t_3\rightarrow t_1$.
Since $G$ is weakly connected, we may assume that $t_1$
has a neighbour in $B$; and by \ref{nbrs} applied taking $v=t_1$, there is an edge $ab$ of $G^-$ where both $a,b$ are adjacent to $t_1$
in $G^-$, and $b\in B$ and $a\notin B$. Consequently $a$ is one of $t_2,t_3$, and so we may assume that
$a=t_2$ (without loss of generality). Since every cycle of $G^-$ of length three is a 
tricycle of $G$ (because $G$ admits a 3-ring), it follows that $t_2b, bt_1\in E(G)$. By \ref{diamond}, $t_3$ has no neighbour
in $B$. From \ref{diamond}, $b$ is the only vertex in $B$ that makes a tricycle with $t_1t_2$. This proves \ref{oneper}.~\bbox

Let us complete the proof of \ref{all3cycled}, which we restate:
\begin{thm}\label{all3cycled3}
Every unbreakable 3-cyclic digraph can be made by special edge sums, starting from 3-pinched digraphs.
\end{thm}
\Proof We show by induction on $|G|$ that if $G$ is an unbreakable 3-cyclic digraph, then 
$G$ can be made by special edge sums, starting from 3-pinched digraphs. From \ref{all3cycled2} we may assume that $G$ is diwheel-free.
If $|G|\le 3$ the result is true, so we assume that $|G|\ge 4$, and hence there is a tricycle  $T$, and a component
$B$ of 
$G^-\setminus V(T)$. By \ref{oneper}, $G^-$ is not 3-connected, and by \ref{findspecial}, the result follows inductively. This proves \ref{all3cycled3}.~\bbox

But if we proceed with more care, we can prove a decomposition result in the diwheel-free case which is more helpful.
\begin{thm}\label{peripheral}
Let $G$ be an unbreakable, 3-cyclic, diwheel-free digraph, with $|G|\ge 4$. 
Then there is an edge $t_1t_2$ of $G$ such that $G^-\setminus \{t_1,t_2\}$ has more than one component, and at most one of them has more than one vertex.
\end{thm}
\Proof
Choose a tricycle $T$, and a component $B$ of $G^-\setminus V(T)$, choosing $T,B$ such that $B$ is maximal. (This is possible since 
$|G|\ge 4$.) Let $T$ be $t_1\rightarrow t_2\rightarrow t_3\rightarrow t_1$. By \ref{oneper}, one of $t_1,t_2,t_3$
has no neighbour in $B$, say $t_3$; and then there is a unique $b\in B$ that makes a tricycle with $t_1t_2$.
For every tricycle $T'$ of $G$ with no vertex in $B$, the choice of $T,B$ implies that
$B$ is a component of $G^-\setminus V(T')$, and so $t_1,t_2\in V(T')$. Since every vertex in $V(G)\setminus (V(T)\cup B)$
belongs to a tricycle, which is necessarily disjoint from $B$, it follows that every vertex in $V(G)\setminus (V(T)\cup B)$
is adjacent from $t_1$ and to $t_2$, and so all such vertices make singleton components of $G^-\setminus \{u,v\}$. 
This proves \ref{peripheral}.~\bbox

We deduce \ref{safeadd}, which we restate:
\begin{thm}\label{safeadd2}
For a digraph $G$, the following are equivalent:
\begin{itemize}
\item $G$ is unbreakable, 3-cyclic and diwheel-free;
\item $G$ is safely buildable.
\end{itemize}
\end{thm}
\Proof
That the second bullet implies the first is easy and we omit it. We prove the first implies the second by induction on $|G|$.
We proceed by induction on $|G|$; so we assume $G$ has at least four vertices, and the theorem holds for all digraphs with fewer 
vertices.  Choose $t_1,t_2$ as in \ref{peripheral}. If $G^-\setminus\{t_1,t_2\}$ has no non-singleton component then $G$ is safely buildable,
so we assume that $B$ is the unique non-singleton component of $G^-\setminus\{t_1,t_2\}$. Let $X = V(G)\setminus (B\cup \{t_1,t_2\})$,
and let $H=G\setminus X$ and choose $t_3\in X$. It follows that $H$ is unbreakable, 3-cyclic and diwheel-free,
and hence safely buildable, from the inductive hypothesis. From \ref{oneper}, there is a unique $b\in B$ that makes a tricycle
with $t_1t_2$. If both $bt_1$ and $t_2b$ belong to a second tricycle, let $b\rightarrow t_1\rightarrow s_1\rightarrow b$ and 
$t_2\rightarrow b\rightarrow s_2\rightarrow t_2$ be tricycles; then $s_1,s_2,t_3$ are all different, and 
$$b\rightarrow s_2\rightarrow t_2\rightarrow t_3\rightarrow t_1\rightarrow s_1\rightarrow b$$
is a directed cycle of length six, a contradiction. So from the symmetry, we may assume that $bt_1$ is in no tricycle of $G$
except the one through $t_2$. Hence, in $H$, the edges $bt_1$ and $t_1t_2$ are both in unique tricycles; and so $t_1$ has degree two
in $H$, by \ref{nbrs}. Consequently $G$ is safely buildable. This proves \ref{safeadd2}.~\bbox

This has some pleasing consequences:
\begin{itemize}
\item If $G$ is safely buildable, then it follows inductively that it admits a 3-ring $(A_1,A_2,A_3)$ in which each of 
$G^-[A_1\cup A_2], G^-[A_2\cup A_3], G^-[A_3\cup A_1]$ is a tree. Since any digraph that admits such a 3-ring is clearly diwheel-free, 
this shows the equivalence of the first two bullets of \ref{nowheel}.
\item If $G$ is safely buildable, then it can certainly be built according to the third bullet of \ref{nowheel}; and since digraphs bulit in this way are diwheel-free (because diwheels have minimum degree three), this shows the equivalence of the first and third bullets of 
\ref{nowheel}.
\item Every unbreakable 3-cyclic digraph $G$ can be built by adding ears, of length one or two, starting 
from a tricycle, by \ref{ears}; and so has at least $3|G|-3$ edges. If $G$ is in addition diwheel-free, then it can be built just by adding ears of length two, and so has exactly $3|G|-3$ edges. Thus an unbreakable, 3-cyclic digraph $G$ is diwheel-free
if and only if it has exactly $3|G|-3$ edges. This shows the equivalence of the fourth bullet of \ref{nowheel} with the others.
\end{itemize}

\section{3-annular digraphs}

Where do the 3-annular digraphs fit in the general picture? One way to answer this is to consider what happens as we build $G$ as in
\ref{safeadd}.
As we add a new vertex $v$, it is made adjacent to some previous vertex $u$ that had degree two; let us call $u$ the ``parent''  of $v$.
This defines a tree $T$ say. $T$ could be any tree (except the vertex $r$
must have degree one), but $G$ is 3-annular if and only if $T$ is a caterpillar  (that is, some path of $T$ with one end $r$ 
contains all vertices with degree at least two in $T$). This can be derived from the equivalence of the first two bullets of 
\ref{activity}
below, but we omit the proof.

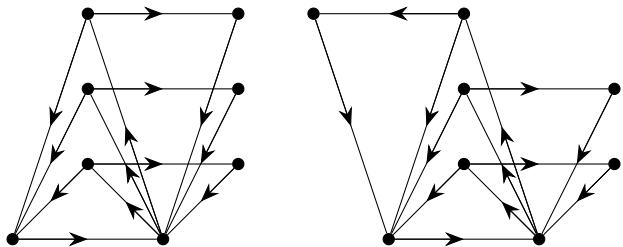
\begin{figure}[ht]
\centering

\begin{tikzpicture}[auto=left]
\tikzstyle{every node}=[inner sep=1.5pt, fill=black,circle,draw]

\node (x) at (0,0) {};
\node (y) at (2,0) {};
\node (a1) at (1,1) {};
\node (a2) at (1,2) {};
\node (a3) at (1,3) {};
\node (b1) at (3,1) {};
\node (b2) at (3, 2) {};
\node (b3) at (3,3) {};

\foreach \from/\to in {x/y,a1/x,a2/x,a3/x,y/a1,y/a2,y/a3, a1/b1,a2/b2,a3/b3, b1/y,b2/y,b3/y}
\draw [] (\from)--(\to);

\def\s{3/2}
\tikzstyle{every node}=[];
\draw [-{Stealth[scale=\s]}] (x) -- (1,0);
\draw [-{Stealth[scale=\s]}] (1,1) -- (.5,.5);
\draw [-{Stealth[scale=\s]}] (1,2) -- (.5,1);
\draw [-{Stealth[scale=\s]}] (1,3) -- (.5,1.5);
\draw [-{Stealth[scale=\s]}] (y) -- (1.5,.5);
\draw [-{Stealth[scale=\s]}] (y) -- (1.5,1);
\draw [-{Stealth[scale=\s]}] (y) -- (1.5,1.5);
\draw [-{Stealth[scale=\s]}] (a1) -- (2,1);
\draw [-{Stealth[scale=\s]}] (a2) -- (2,2);
\draw [-{Stealth[scale=\s]}] (a3) -- (2,3);
\draw [-{Stealth[scale=\s]}] (b1) -- (2.5,.5);
\draw [-{Stealth[scale=\s]}] (b2) -- (2.5,1);
\draw [-{Stealth[scale=\s]}] (b3) -- (2.5,1.5);

\tikzstyle{every node}=[inner sep=1.5pt, fill=black,circle,draw]
\node (x) at (5,0) {};
\node (y) at (7,0) {};
\node (a1) at (6,1) {};
\node (a2) at (6,2) {};
\node (a3) at (6,3) {};
\node (b1) at (8,1) {};
\node (b2) at (8, 2) {};
\node (b3) at (4,3) {};

\tikzstyle{every node}=[];
\foreach \from/\to in {x/y,a1/x,a2/x,a3/x,y/a1,y/a2,y/a3, a1/b1,a2/b2,a3/b3, b1/y,b2/y,x/b3}
\draw [] (\from)--(\to);

\draw [-{Stealth[scale=\s]}] (x) -- (6,0);
\draw [-{Stealth[scale=\s]}] (6,1) -- (5.5,.5);
\draw [-{Stealth[scale=\s]}] (6,2) -- (5.5,1);
\draw [-{Stealth[scale=\s]}] (6,3) -- (5.5,1.5);
\draw [-{Stealth[scale=\s]}] (y) -- (6.5,.5);
\draw [-{Stealth[scale=\s]}] (y) -- (6.5,1);
\draw [-{Stealth[scale=\s]}] (y) -- (6.5,1.5);
\draw [-{Stealth[scale=\s]}] (a1) -- (7,1);
\draw [-{Stealth[scale=\s]}] (a2) -- (7,2);
\draw [-{Stealth[scale=\s]}] (b3) -- (5,3);
\draw [-{Stealth[scale=\s]}] (b1) -- (7.5,.5);
\draw [-{Stealth[scale=\s]}] (b2) -- (7.5,1);
\draw [-{Stealth[scale=\s]}] (x) -- (4.5,1.5);

\end{tikzpicture}

\caption{Two of the branchers } \label{fig:branchers}

\end{figure}

Finally, we would like to give a characterizations of the 3-annular digraphs. The {\em branchers} are the two digraphs of Figure \ref{fig:branchers},
together with the two digraphs obtained from these by reversing the direction of all edges.
If $uv$ is an edge of a digraph $G$, its {\em activity} is the number of components of $G^-\{u,v\}$ that have more than one vertex. 
We have the following, which implies \ref{planar}:
\begin{thm}\label{activity}
Let $G$ be an unbreakable diwheel. Then the following are equivalent:
\begin{itemize}
\item $G$ is 3-annular;
\item  $G$ is 3-cyclic and diwheel-free, and every edge has activity at most two;
\item $G$ is 3-cyclic and diwheel-free, and contains none of the four branchers.
\end{itemize}
\end{thm}
\Proof
The four branchers are not 3-annular, so the first bullet implies the third. Next, we show that 
the third implies the second. To see this, observe that if
$B$ is a non-singleton component $G^-\setminus \{u,v\}$ for some edge $uv$, then there is a tricycle consisting of $u,v$ and a vertex $b\in B$; and since $|B|>1$, $b$ has a neighbour in $B$, and \ref{nbrs} implies that there is a tricycle containing $b$, one of $u,v$, and some other
vertex of $B$. This gives two cases, depending whether $ub$ or $vb$ is in the second tricycle. If $uv$ has activity at least 
three, then the same holds for three choices of $B$, and so $G$ contains one of the four branchers.

It remains to show that the second bullet implies the first, and we prove this by induction on $|G|$. Let $G$ be unbreakable, 
3-cyclic and diwheel-free, such that every edge has activity at most two. We may assume that $|G|\ge 4$.
By \ref{peripheral}, we can choose an edge $t_1t_2$ of $G$ such that $G^-\setminus \{t_1,t_2\}$ has more than one component, and at 
most one of them is not a singleton. If all of them are singletons, the result holds, so we assume that $B$ is the 
unique non-singleton component of $G^-\setminus \{t_1,t_2\}$. Choose $t_3\in X$.  As in the proof of \ref{safeadd2},
there is a unique $b\in B$ making a tricycle with $t_1t_2$; and we may assume that $t_1$ has degree two in $G\setminus X$.
Moreover, $H=G\setminus X$ is unbreakable, 3-cyclic and diwheel-free, and so admits a 3-annular drawing, from the inductive hypothesis. 
From the hypothesis, the edge $bt_1$ has activity at most two, and since $\{t_2\}\cup X$ is the vertex set of a 
non-singleton component of $G^-\setminus \{b,t_1\}$, it follows that $H\setminus \{t_1,t_2,b\}$ has at most non-singleton component. Since $H$ admits
a 3-annular drawing, and $H\setminus \{t_1,t_2,b\}$ has at most non-singleton component, it follows that we can rearrange the drawing 
so that $t_1,t_2,b$ are the vertices of the outermost tricycle. But then we can add $X$ to the drawing, and make a 3-annular drawing of $G$.
This proves \ref{activity}.~\bbox

\section{More about weightable digraphs}

In this section we prove the claims made in the introduction about weightable digraphs. If $G$ is a digraph, a {\em weighting} is a 
map $w:E(G)\to \mathbb{R}$ such that $\sum_{e\in E(C)} w(e) = 1$ for every directed cycle $C$ ($\mathbb{R}$ is the set of real numbers.)
Thus $G$ is weightable if and only if it admits a weighting.
We begin with:
\begin{thm}\label{realtoint}
If $G$ is weightable, then it admits a weighting that is integer-valued.
\end{thm}
\Proof
We may assume by induction on $|G|$ that $G$ is unbreakable. Choose a maximal unbreakable subdigraph $H$ such that there is a weighting $w$ of $G$ satisfying $w(e)\in \mathbb{Z}$ for each $e\in E(H)$. ($\mathbb{Z}$ is the set of integers.)  
Suppose that $H\ne G$; then, by \ref{ears},  there is an ear $P$ in the sense of \ref{ears}, from $u$ to $v$ say. Since $H$ is unbreakable,
there is a directed path $Q$ of $H$ from $v$ to $u$. Let $C$ be the directed cycle $P\cup Q$, and let the vertices of $P$ be 
$u=p_1\CC p_n=v$ in order. Let $e_i$ be the edge $p_ip_{i+1}$ for $1\le i<n$, and for $1\le i\le n$ let $P_i$ be the subpath of $P$ from $u=p_1$ to $p_i$. Choose $i\in \{1\LL n\}$ 
maximal such that there is a weighting $w$ of $G$ satisfying that $w(e)\in \mathbb{Z}$ for each $e\in E(H)\cup E(P_i)$ (this is 
possible since we may take $i = 0$ if necessary). By the maximality of $H$ it follows that $i\le n-1$; and if $i=n-1$, then since
$w(C)=1$ and $w(e)\in \mathbb{Z}$ for each edge $e\in E(C)$ except $e_{n-1}$, it follows that $w(e_{n-1})\in \mathbb{Z}$,
a contradiction to the definition of $i$. So $i\le n-2$. Define $w':E(G)\to \mathbb{R}$ by: $w'(e)=w(e)$ for each $e\in E(G)$ not 
incident with $p_{i+1}$; $w'(e)=w(e)-w(e_i)$ for each edge $e$ 
with head $p_{i+1}$; and $w'(e) = w(e)+w(e_i)$ for each edge with tail $p_{i+1}$. Then $w'$ is a weighting, and 
$w'(e)\in \mathbb{Z}$ for each $e\in E(H)\cup E(P_{i+1})$, contradicting the choice of $i$. This proves that $H=G$, and so proves \ref{realtoint}.~\bbox

\begin{thm}\label{inttoone}
If $G$ is weightable, it admits a weighting that is $\{0,1\}$-valued.
\end{thm}
\Proof
We may assume that every edge is in a directed cycle, since we may delete any edges not in directed cycles.
By \ref{realtoint}, there is an integer-valued weighting $w$. Define $\phi(w)$ to be the sum of $w(e)$ over all $e\in E(G)$ with 
$w(e)<0$; and choose $w$ to maximize $\phi(w)$.
Suppose that $\phi(w)<0$, and choose $uv\in E(G)$ with $w(uv)<0$. Let $X$ be the set of all vertices $x$ such that there is
a directed path $P$ of $G$ from $v$ to $x$ where $w(e)\le 0$ for all edges $e\in E(P)$. Thus $v\in X$, and $u\notin X$, since adding $uv$ to such a path
from $v$ to $u$ would give a directed cycle with negative total weight. Let $A$ be the set of all edges with tail in $X$ and head in 
$V(G)\setminus X$, and $B$ the set of edges with head in $X$ and tail in $V(G)\setminus X$. So $uv\in B$; and $w(e)\ge 1$ for each $e\in A$,
from the definition of $X$ and since $w$ is integer-valued. Define $w':E(G)\to \mathbb{R}$ by: $w'(e)=w(e)$ for each 
$e\in E(G)\setminus (A\cup B)$; $w'(e) = w(e) -1$ for each $e\in A$; and $w'(e) = w(e) + 1$ for each $e\in B$. Then $w'$ is an 
integer-valued weighting, and $\phi(w')>\phi(w)$, a contradiction. This proves that $\phi(w)=0$, and so $w(e)\ge 0$ for every edge $e$.
Since every edge is in a directed cycle, and $w$ is a weighting, it follows that $w(e)\in \{0,1\}$ for each edge $e$. This proves \ref{inttoone}.~\bbox

Next, we will characterize the minimal digraphs that are not weightable. The proof is almost exactly the same as a proof in 
\cite{evendigraph}, except that there we were characterizing the minimal digraphs that were not weightable modulo two; that is, those that did not admit a map $w:E(G)\to \{0,1\}$ such that $w(C)$ is odd for every directed cycle $C$. That proof is quite long, and great swathes of it
need minimal or no changes, so we will just indicate the changes needed in the proof of \cite{evendigraph}.

Let $k\ge 2$, and let $G$ be the digraph with vertex set $\{v_1\LL v_k\}$ and edges $v_iv_{i+1}$ and $v_{i+1}v_i$ for $1\le i\le k$, reading subscripts modulo $k$. We call $G$ a {\em $k$-double-cycle}. (A 2-double-cycle is just a directed cycle of length two.)
A {\em weak $k$-double-cycle} is a digraph formed by the union of $k$ directed cycles $C_1\LL C_k$ arranged in circular order, where
each vertex belongs to at most two of $C_1\LL C_k$, and $C_i\cap C_{i+1}$ is a directed path for $1\le i\le k$ (reading subscripts 
modulo $k$) and $C_i$ is vertex-disjoint from $C_j$
if $j\ne i-1,i,i+1$ modulo $k$. In \cite{evendigraph} we showed that a digraph is weightable modulo two if and only if no subdigraph 
is a  weak $k$-double-cycle where $k\ge 3$ and odd. Here we will show that:
\begin{thm}\label{doublecycle}
A digraph $G$ is weightable if and only if no subdigraph is a weak $k$-double-cycle with $k\ge 3$.
\end{thm}

If $C$ is a directed cycle, let ${\bf C}$ be its characteristic vector in $\mathbb{R}^{E(G)}$; so ${\bf C}_e = 1$ if 
$e\in E(C)$ and ${\bf C}_e = 0$ otherwise. We view ${\bf C}$ as a column vector.
Let us say a {\em real directed cycle basis} of a digraph $G$ is a set $\mathcal{C}$ of directed cycles such that:
\begin{itemize}
\item for every directed cycle $D$ of $G$, ${\bf D}$ is a linear combination of the vectors ${\bf C}\;(C\in \mathcal{C})$; and
\item the vectors ${\bf C}\;(C\in \mathcal{C})$ are linear independent.
\end{itemize}

The proof of Proposition (2.1) in~\cite{evendigraph} shows without any change 
that:
\begin{thm} \label{new2.1}
Every strongly connected  digraph $G$ has a real directed cycle basis $\mathcal{C}$.
Moreover, if $A\subseteq E(G)$ is such that $G\setminus A$ is strongly connected, then $\mathcal{C}$ can be chosen
such that each edge of $A$ is in precisely one cycle of $\mathcal{C}$.
\end{thm}

Moreover, the later part of the proof of theorem (4.1) in~\cite{evendigraph} (the section from step (4) onwards) shows the following
(although this is not stated explicitly in that paper):
\begin{thm}\label{getdoublecycle}
Let $G$ be a strongly 2-connected digraph. If for every $F\subseteq E(G)$ such that $G\setminus F$ is strongly connected, there is a 
directed cycle containing $F$, then $G$ is a $k$-double-cycle for some $k\ge 2$. 
\end{thm}

\noindent{\bf Proof of \ref{doublecycle}.\ \ }
We mimic the proof of theorem (4.1) in~\cite{evendigraph}, which we call the ``old proof''. 
Suppose that $G$ is a minimal counterexample. As in the old proof, it follows that $G$ is strongly connected and $|G|\ge 4$. 
The argument in step (2) of the old proof also works (except the weightings in the proof are in the sense of this paper),
and we deduce that $G$ is strongly 2-connected. Step (3) of the old proof needs some adjustment. We need to show:
\\
\\
(1) {\em If $F\subseteq E(G)$ is such that $G\setminus F$ is strongly connected, then there is a directed cycle of $G$ containing $F$.}
\\
\\
By \ref{new2.1}, there is a real directed cycle basis $\mathcal{C}$ such that every edge in $F$ belongs to exactly one cycle of 
$\mathcal{C}$. Let ${\bf M}$ be the real matrix with columns indexed by $\mathcal{C}$ and where the column indexed by $C\in \mathcal{C}$ is 
the vector ${\bf C}$. Since the columns of ${\bf M}$ are linearly independent, there is a vector ${\bf w}\in \mathcal{R}^{E(G)}$ such that 
${\bf M}^{T}{\bf w}={\bf 1}_{\mathcal{C}}$, where superscript $T$ means transpose, and ${\bf 1}_{\mathcal{C}}$ is the all-ones vector in $\mathcal{R}^{\mathcal{C}}$. Since $G$ is not weightable, there is 
a directed cycle $D$ such that ${\bf D}^T{\bf w}\ne 1$. Since $\mathcal{C}$ is a real directed  cycle basis, there is a 
vector
${\bf f}\in \mathbb{R}^{\mathcal{C}}$ such that ${\bf M}{\bf f} = {\bf D}$. Let $\mathcal{C}'$ be the set of cycles $C\in \mathcal{C}$ 
with ${\bf f}_C\ne 0$; let $G'$ be the union of the cycles in $\mathcal{C'}$ (and consequently $D$ is a cycle of $G'$); and let ${\bf N}$ be the submatrix of ${\bf M}$ formed by the columns
${\bf C}$ with $C\in \mathcal{C}'$. Thus ${\bf N}{\bf f} = {\bf D}$, and ${\bf N}^{T}{\bf w}={\bf 1}_{\mathcal{C}'}$, where 
${\bf 1}_{\mathcal{C}'}$ is the all-ones vector in $\mathbb{R}^{\mathcal{C}'}$.
Suppose that $G'\ne G$, and so $G'$ is weightable, from the minimality of $G$. Hence 
there is a vector 
${\bf w}'\in \mathcal{R}^{E(G)}$ such that ${\bf N}^{T}{\bf w'}={\bf 1}_{\mathcal{C}'}$, and ${\bf D}^T{\bf w}'= 1$.
Consequently ${\bf N}^{T}({\bf w}'-{\bf w})={\bf 0}$, and ${\bf D}^T({\bf w}'-{\bf w})\ne 0$. But this contradicts that 
${\bf N}{\bf f} = {\bf D}$, and so proves that $G'=G$. Let $e\in F$, and let $C$ be the unique member of $\mathcal{C}$ that contains $e$.
Since $e\in E(G')$, it follows that $C\in \mathcal{C}'$, and so $f_C\ne 0$. Since ${\bf M}{\bf f} = {\bf D}$, it follows that $e\in E(D)$.
This proves (1).

\bigskip
From (1) and \ref{getdoublecycle} it follows that $G$ is a $k$-double-cycle for some $k\ge 2$. Since the 2-double-cycle is weightable, we deduce that $k\ge 3$. This proves \ref{doublecycle}.~\bbox

\section*{Acknowledgements}
Thanks to a referee for a very careful reading and many helpful comments.
Thanks also to Stephen Bartell for interesting discussions on extensions of this problem; and to Alex Scott and Tung Nguyen 
for discussions on \ref{inttoone}.

\end{document}